\documentclass[12pt,a4paper]{article}

\usepackage{amsmath,amsthm,amssymb,amscd,enumerate}
\usepackage[english]{babel}
\usepackage[latin1]{inputenc}
\usepackage{tikz}
\usepackage{listings}
\usepackage{graphicx,graphics}
\usepackage{bbm}
\usepackage{tocbibind}
\usepackage{fancyhdr}
\usepackage{pgf}

\begin{document}

\title{
An ecoepidemic food chain with the disease at the intermediate trophic level
}
\author{Lorenzo Bramante, Simone Maiolo, Ezio Venturino\\
Dipartimento di Matematica ``Giuseppe Peano'',\\
Universit\`a di Torino, \\via Carlo Alberto 10, 10123 Torino, Italy,\\
email: ezio.venturino@unito.it}
\date{}

\maketitle

\begin{abstract}
We consider a three-level food chain in which an epidemics
affects the intermediate population. Two models are presented, respectively either allowing
for unlimited food supply for the bottom prey, or instead assuming for it a logistic growth.
Counterintuitive results related to the paradox of enrichment are obtained,
showing that by providing large amounts of food to the bottom prey, the top
predator and the disease in suitable situations can be eradicated.
\end{abstract}


\section{Introduction}

Food chains constitute a very common ecological situation.
For an earlier model of this kind, see for instance \cite{[9]}.
In \cite{Selakovic}, a whole wealth of real
life examples are presented and discussed.
In particular,
\cite{Du} contains  a description of a cascade of diseases that moved from rinderpest for cattle to wild animals
and then because of the death of these herds, caused also human diseases including smallpox.
The study focuses on the reforestation of the Serengeti Woodlands along the past century.
In it not only the role of fires in altering the landscape is described, but also the more relevant one of elephants.
These animals stip the bark of trees and break their branches, contributing substantially to their decline.
This behavior is common to several herbivores at all latitudes, \cite{Tamburino1,Tamburino2}.
Also, the influence of the tsetse fly as carriers of the related infections of trypanosomes are highlighted.
This disease affects large animals like cattle, but does not harm the small herbivores. It is the cause in man
of the ``sleeping sickness'' disease. Overgrazing of cattle removes high quantities of grass and makes
fires occurrence less frequent so that bushes can regrow.
There seems to be a cyclic behavior among these phases (trees and tsetse, grass and fires) along the past century,
evidence of a dynamic ecosystem, very much intricated.

In another investigation about the Serengeti Woodlands, \cite{HS},
it is observed again that the reforestation
is tightly related to diseases of wild animals
populating or invading
that environment, and therefore diseases, in this specific case
rinderpest, play an essential role in regulating the ecosystem. Occurrence of epidemics among the
herbivores has far reaching consequences not just for the animals, but also inflence the whole ecosystem via a kind of chain reaction.
At the same time wild fires clearly control the canopy, changing the size of Carbon stored in the soil and the biomass.

Mathematical models for diseases affecting interacting populations are known since two decades
at least, \cite{HF,BC,V94,V95}, and involve interactions of every possible kind,
\cite{CA99,V01,V02,V07,SH} and various other modeling assumptions, \cite{AAMC,HV07}.
Ecoepidemiology, see Chapter 7 of \cite{MPV}, is the study of such ecosystems.
So far, investigations have confined themselves essentially to simple systems, mainly two
intermingling populations with one disease affecting one of them. But very recently epidemics
in food chains have been considered, \cite{ECOCOM}.

In this paper we continue the investigations of \cite{ECOCOM},
in which the epidemics propagates instead at the lowest trophic level,
by considering the infected individuals to be predators
on the bottom prey,
but also subject themselves to being hunted by a top predator.

Two models are here presented, after that the underlying basic demographic system is analysed,
and then in turn the Malthus and the logistic versions of the ecoepidemic food chains
are studied. A final interpretation of the results concludes the paper.

\section{The general model}

We consider a three trophic level food chain, composed by the populations $P,$ $H$ and $V$, in which the intermediate population is subject to
a disease transmissible by contact at rate $\beta$. We therefore partition it into the two sets of susceptibles $S$ and infected $I$. We assume the disease
to be unrecoverable. Also, it is confined to the population $H$ and cannot be trasmitted either to its predators $P$ or its prey $V$.
The infected are weakened by the disease so much so as to be unable to
extert any pressure on the population $V$, nor to feel any such pressure
from the healthy individuals of their own population; they can be captured by
the top predators but do not cause them any harm.
The top predators do not have any food sources other than their prey $H$.

The model in the logistic formulation is
\begin{eqnarray}\label{eq:sistemasemplificato}
\begin{split}
\frac{dP}{dt}&=P \left( gI+fS-\tau \right),\\
\frac{dS}{dt}&=S \left( lV-\beta I-qP-\mu \right),\\
\frac{dI}{dt}&=I \left(\beta S-cP-\nu \right),\\
\frac{dV}{dt}&=V\left[ r\left(1-\frac{V}{K}\right)-bS\right].
\end{split}
\end{eqnarray}
The first equation states that the top predators in absence of $H$ would die out at an exponential rate. They can survive by predation on the next
trophic level, and as stated are not harmed by eating infected individuals. In the second equation we find the dynamics of the healthy individuals
of the intermediate population. They reproduce as long as they can feed on the lower population $V$, and leave this class either by becoming infected,
or by mortality, whether be it natural or induced by their capture from the top predators.
The next equation contains the infected behavior; the only input is due to the $S$ individuals that become diseased upon ``successful'' contact
with a disease-carrier. Infected leave this class if they are hunted by the $P$'s, or by mortality, which can also be induced by the disease.
The last equation states that the lower population in the trophic level reproduces logistically and is hunted only by the healthy individuals
of the upper trophic level.

The meaning of the parameters is as follows:
$r$ denotes the reproduction rate of the population $V$ and $K$ is its respective carrying capacity;
$b$ is the hunting rate of $S$ on $V$;
$\beta$ is the disease incidence rate;
$c$ and $q$ are the predation rate of $P$ on $I$ and $S$ respectively;
$\mu$ is the population $H$ natural mortality rate, while $\nu:=\mu+\mu_{0}$ represents the mortality rate for the infected, which includes
the disease-related mortality $\mu_{0}$; finally $\tau$ is the mortality rate for the population $P$.
In view that not all prey are converted into predator's biomass, we have the restrictions
\begin{equation}\label{param}
g<c, \quad f<q, \quad l<b. 
\end{equation}

When the resources for $V$ are unlimited, i.e. for $K\rightarrow \infty$, we have the Malthus case, i.e. (\ref{eq:sistemasemplificato}) simplifies as follows:
\begin{eqnarray}\label{sistMalthus}
\begin{split}
\frac{dP}{dt}&=gPI+fPS-\tau P,\\
\frac{dS}{dt}&=lSV-\beta SI-qSP-\mu S,\\
\frac{dI}{dt}&=\beta IS-cIP-\nu I,\\
\frac{dV}{dt}&=rV-bVS.
\end{split}
\end{eqnarray}

The Jacobian $J$ of (\ref{eq:sistemasemplificato}) is
\begin{equation}\label{Jac}
\left[
\begin{array}{cccc}
gI+fS-\tau & fP & gP & 0\\
-qS & lV-\mu-\beta I-qP & -\beta S & lS \\
-Ic & \beta I & -\nu+\beta S-cP & 0 \\
0 & -Vb & 0 & r(1-\frac{V}{K})-bS-\frac{Vr}{K}
\end{array}
\right]
\end{equation}

Note that the Jacobian of (\ref{sistMalthus}) contains a modification only
in the last term of the last equation, namely
\begin{equation}\label{eq:jacobMalthus}
J=
\left[
\begin{array}{cccc}
gI+fS-\tau & fP & gP & 0\\
-qS & lV-\mu-\beta I-qP & -\beta S & lS \\
-Ic & \beta I & -\nu+\beta S-cP & 0 \\
0 & -Vb & 0 & r-bS
\end{array}
\right].
\end{equation}

\section{The disease-free model}
We replace the two intermediate equations of (\ref{eq:sistemasemplificato}) by their total population
$Q=S+I$, and observing that there are no infected in this case, in fact
$Q=S$, thus obtaining the equation 
$$
\frac{dQ}{dt}=- \mu Q - qQP +l QV
$$
Also, the Jacobian becomes a $3\times 3$ matrix. Corresponding changes occur in
(\ref{sistMalthus}) and (\ref{eq:jacobMalthus}).

The system has only three meaningful equilibria, since the origin is unconditionally unstable.
The bottom prey-only
equilibrium $D_1=(0,0,K)$ exists only in the logistic case. The top-predator-free equilibrium
$\widehat D=\left( 0, \widehat Q, \widehat V\right)$,
\begin{eqnarray*}
\widehat Q=\frac rb\left( 1-\frac {\mu}{lK}\right), \quad \widehat V=\frac {\mu}l
\end{eqnarray*}
and the coexistence equilibrium
$D^*=\left( P^*, Q^*, V^*\right)$, whose population values are
\begin{eqnarray*}
P^*=\frac 1q \left[ l K \left( 1- \frac {b\tau}{rf}\right) - \mu \right], \quad
Q^*=\frac {\tau}f, \quad
V^*=K \left( 1- \frac{b\tau}{rf}\right).
\end{eqnarray*}

Now, $D_1$ is stable if
\begin{equation}\label{D1_stab}
1>\frac {lK}{\mu} \equiv \rho_1,
\end{equation}
while $\widehat D$ is feasible in the opposite case,
\begin{equation}\label{Dhat_feas}
\rho_1 \ge 1.
\end{equation}
Thus we have a transcritical bifurcation. $\widehat D$ is stable for
\begin{equation}\label{Dhat_stab}
1 >  \frac {fr}{b \tau}\left( 1-\frac {\mu}{lK}\right) \equiv \rho_2.
\end{equation}
The opposite condition provides instead feasibility for $D^*$:
\begin{equation}\label{D*_feas}
\rho_2\ge 1,
\end{equation}
thus we have another transcritical bifurcation. Stability of $D^*$ holds unconditionally,
whenever the equilibrium is feasible, as the Routh-Hurwitz conditions become
\begin{equation}\label{D*_stab}
\frac rK V^*>0, \quad
fq\frac rK P^*Q^*V^*>0, \quad
bl V^*>0.
\end{equation}

\section{The Malthus case}
For the ecoepidemic cases, we analyse at first the particular case of (\ref{eq:sistemasemplificato}).

The possible equilibria are the following points: since the system is homogeneous, the origin trivially satisfies it, $\widetilde E_0=(0,0,0,0)$. Then we have
$$
\widetilde E_1=\left( 0,\frac{r}{b},0,\frac{\mu}{l}\right)
$$
which is always feasible.
Finally, coexistence is obtained at the level
\begin{equation}
\widetilde E^*=\left( \frac{\beta r-\nu b}{bc},\frac{r}{b},\frac{b\tau-rf}{gb},
\frac{\beta rgq-\beta rfc+bg\mu c-bgq\nu+b\tau\beta c}{gbcl}\right).
\end{equation}
Feasibility implies that all the following conditions hold
\begin{equation}\label{feas_E16}
\frac {r\beta}{b\nu} \ge 1; \quad \frac {rf}{b\tau} \le 1; \quad
b\frac{\tau\beta c+g\left(\mu c-q\nu\right)}{r\beta\left(fc-gq\right)} \le 1.
\end{equation}

\vspace{0.3cm}
{\textbf{Remark 1}}. It is interesting to note that the bottom population-free point
$$
\left( -\frac{\nu f-g\mu-\tau\beta}{-gq+fc},\frac{g\mu c-gq\nu+\tau\beta c}{\beta(-gq+fc)},
-\frac{\mu cf-q\nu f+q\tau\beta}{\beta\left(-gq+fc\right)},0\right)
$$
intuitively cannot be an equilibrium, since what is called the primary producer, $V$, is wiped out,
and therefore the first trophic level, i.e. populations $S$ and $I$ that feed on it, cannot
thrive any longer, and in turn also the top predator must die out, since the intermediate
population is depleted. This observation has its counterpart in the mathematics, since this point
is not feasible. In fact requiring all its populations to be nonnegative leads to the two mutually exclusive conditions
$$
\frac{g}{c}> -\frac{\tau \beta}{\mu c-q\nu}; \qquad  \frac{g}{c}<-\frac{\tau \beta}{\mu c-q\nu}.
$$
\vspace{0.3cm}

{\textbf{Remark 2}}. Some similar considerations can be made in a few other cases.
In particular note that for the top predator-free subsystem cannot settle to an
equilibrium, quite unexpectedly, because removing the $P$ population and its
related differential equation, we find from the last two equilibrium equations
that $S$ attains the values
$$
S=\frac {\nu}{\beta}, \quad S=\frac rb
$$
which cannot be equal except for a very restrictive condition on the parameters,
that in general does not hold.
\vspace{0.3cm}

Easily, the eigenvalues of the Jacobian (\ref{eq:jacobMalthus}) at the origin are
$-\nu$, $-\tau$, $-\mu$, $r$, from which instability of $\widetilde E_0$ follows.

At $\widetilde E_1$ we find
$$
\pm i\sqrt{\mu r}, \quad \frac{\beta r-\nu b}{b}, \quad \frac{rf-b\tau}{b}.
$$
It follows that if we require
\begin{equation}\label{stab_E7M}
\frac {r\beta}{b\nu} < 1, \quad \frac {rf}{b\tau} <1
\end{equation}
we obtain the neutral, or center, stability, since the remaining two eigenvalues are pure imaginary.
This feature is of course inherited by the fact that the underlying demographic model in this
case is the classical Lotka-Volterra predator-prey model.
On comparing the first stability condition (\ref{stab_E7M}) with the first feasibility
condition (\ref{feas_E16}) of $\widetilde E^*$, we discover a transcritical bifurcation for which
coexistence originates from the equilibrium $\widetilde E_1$ when the latter becomes unstable.

The equilibrium $\widetilde E^*$ exhibits a fourth degree characteristic polynomial,
\begin{equation}\label{char_pol}
\sum _{k=0}^4 a_{4-k}\lambda ^k , \quad a_0=1,
\end{equation}
with known but rather complicated coefficients, which we omit. In any case, we find that $a_1=0$,
so that the very first Routh-Hurwitz condition, $a_1>0$ is not satisfied.
We conclude then that $\widetilde E^*$ is always unstable.

\begin{figure}[ht]
\includegraphics[scale=0.4]{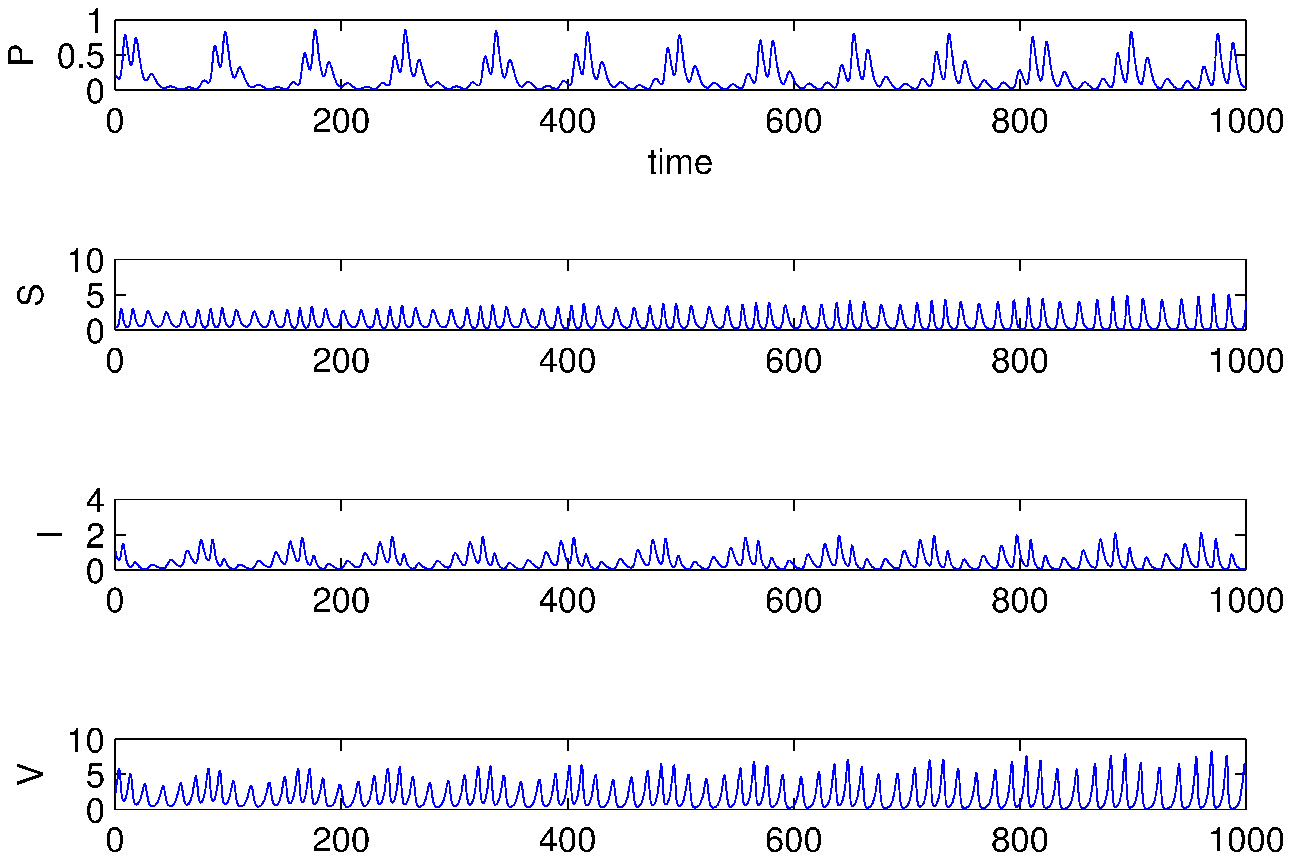}
\includegraphics[scale=0.3]{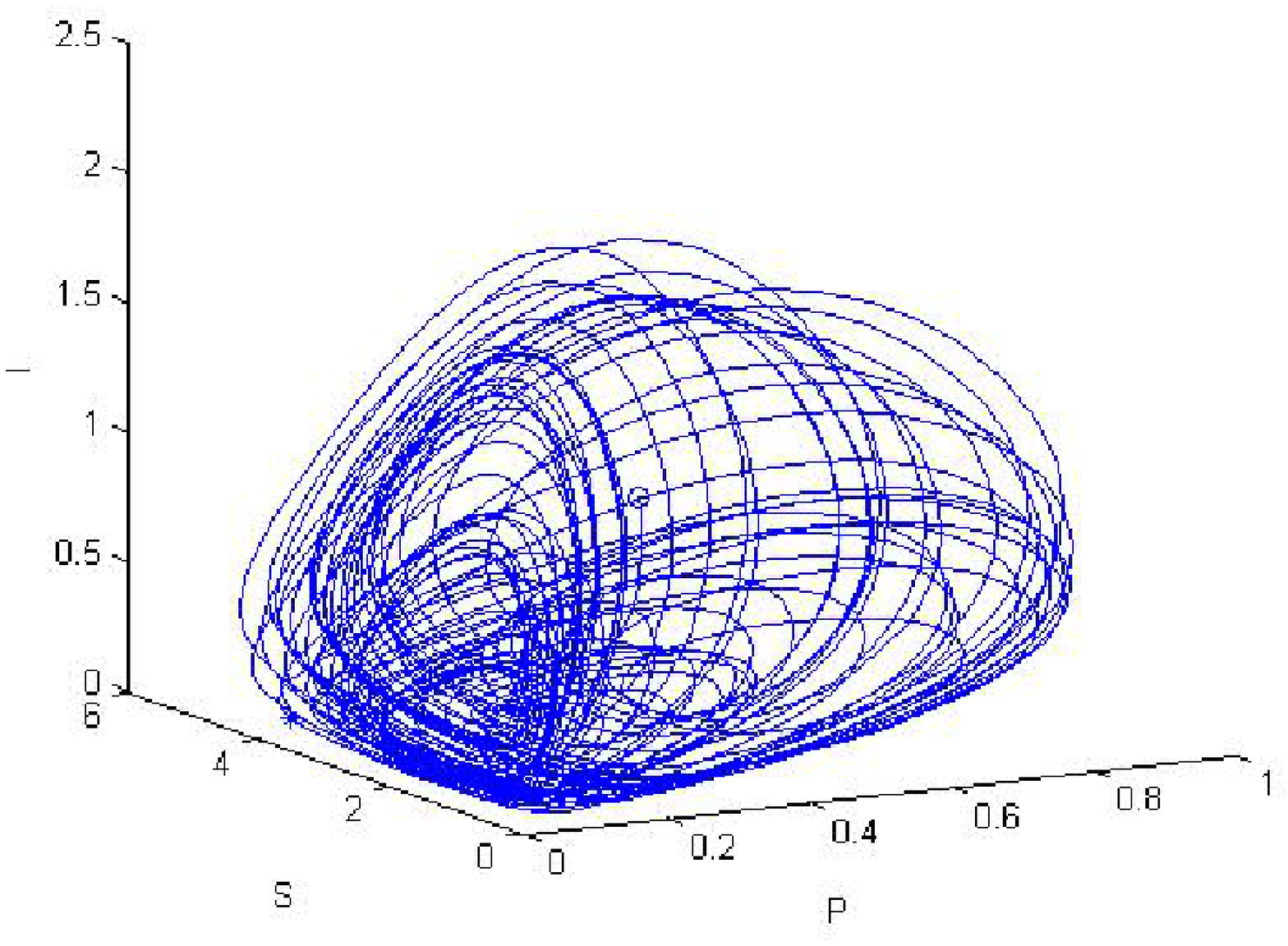}
\caption{Coexistence is attained through persistent oscillations. Left: the populations
as function of time, $P$, $S$, $I$, $V$ from top to bottom.
Right: the $PSI$-subspace phase portrait of the system trajectories.
Parameter values:
$g=0.3$, $f=0.2$, $c=0.4$, $l=0.2$, $q=0.3$, $b=0.4$, $\beta=0.3$, $\tau=0.4$,
$\nu=0.3$, $\mu=0.2$, $r=0.5$.
}
\label{fig:coex_M}
\end{figure}

Coexistence then can only occur at unstable level, i.e. via oscillations. This is
shown in Figure \ref{fig:coex_M} for the parameter values
$g=0.3$, $f=0.2$, $c=0.4$, $l=0.2$, $q=0.3$, $b=0.4$, $\beta=0.3$, $\tau=0.4$,
$\nu=0.3$, $\mu=0.2$, $r=0.5$.

\section{The logistic case}
We now consider (\ref{eq:sistemasemplificato}).
In this case it is possible to show boundedness of the system.

\vspace{0.3cm}
{\textbf{Theorem}}. The system's trajectories are bounded.

\vspace{0.3cm}
{\textbf{Proof}}.
Let us consider the total environment population, $W=P+S+I+V$. Upon summation of
the equations in (\ref{eq:sistemasemplificato}) we obtain, 
$$
\frac {dW}{dt}=(g-c)PI+(f-q)PS+(l-b)SV-\tau P-\mu S -\nu I+rV \left(1 - \frac VK\right).
$$
Recalling the relationships between parameters (\ref{param}), introducing an arbitrary $\theta>0$,
we find
$$
\frac {dW}{dt} + \theta W \le -(\tau -\theta )P-(\mu -\theta)S -(\nu -\theta) I+(r +\theta)V
- \frac rK V^2.
$$
Taking $\theta\le \min \{ \tau, \mu, \nu\}$, the first terms in the above inequality
can be dropped. The last two terms are the parabola $\Psi(V)=V[(r +\theta)-rVK^{-1}]$, whose vertex
lies at the point $(V_0,\Psi^*)=\left( (r +\theta)K(2r)^{-1},(r +\theta)^2K(4r)^{-1}\right)$. It therefore follows
$$
\frac {dW}{dt} + \theta W \le \Psi^*
$$
and upon integration of the corresponding differential equation
we find $W(t)=\Psi^* \theta ^{-1} [1-\exp(-\theta t)] + W(0) \exp(-\theta t)$ so that ultimately
$$
W(t) \le \max \left\{ W(0), \frac {\Psi^*}{\theta} \right\},
$$
as desired.

\vspace{0.3cm}

The equilibria are once again the origin $E_0$
and the equilibrium of the
lowest two trophic levels predator-prey disease-free
subsystem, $E_1$, for which now the susceptible
population level is lower than in the Malthus case $\widetilde E_1$, namely
$$
S_1=\frac rb \left( 1- \frac {\mu}{lK} \right), \quad V_1=\frac{\mu}{l},
$$
and $E^*$, whose components cannot in this case be explicitly evaluated.
Feasibility for $E_1$ holds if (\ref{Dhat_feas})
is satisfied.
In addition, we find the bottom-prey-only equilibrium $E_2=(0,0,0,K)$ and two more points,
$$
E_{3}=\left( 0,\frac{\nu}{\beta},\frac{lKr\beta-lKbv-\mu r\beta}{r\beta^{2}},K\frac{r\beta-b\nu}{r\beta}\right)
$$
the top-predator free equilibrium
and the disease-free equilibrium
$$
E_{4}=\left( \frac{rflK-rf\mu -\tau lKb}{rfq},\frac{\tau}{f},0,K\frac{rf-b\tau}{rf}\right).
$$

To have nonnegative populations, we must impose both the following feasibility conditions
\begin{equation}\label{feas_E12}
r\beta (KL-\mu)> bKl\nu.
\end{equation}
Similarly, for feasibility of $E_{4}$ the parameters must instead satisfy
\begin{equation}\label{feas_E14}
rf (KL-\mu)> bKl\tau.
\end{equation}

The eigenvalues at $E_0$ coincide with those of $\widetilde E_0$,
so that the origin retains its unstable character.

At $E_2$ we find the eigenvalues
$-r$, $-\nu$, $-\tau$, $lK-\mu$. It is stable if (\ref{D1_stab}) holds,
which compared with the feasibility condition for $E_1$, (\ref{Dhat_feas}),
shows again
the existence of a transcritical bifurcation, inherited from the demographic model.

For $E_1$ things change a bit, with respect to $\widetilde E_1$.
Namely the coefficients of the characteristic equation (\ref{char_pol}) are now
\begin{eqnarray*}
a_4=\frac{(-rflK+rf\mu+\tau lKb)r(lK-\mu)\mu(-lKr\beta+lKb\nu+\mu r \beta)}{l^{3}K^{3}b^{2}},\\
a_3=\frac{\mu r}{l^{3}K^{3}b^{2}}(-rfl^{3}K^{3}b+2rfl^{2}K^{2}b\mu+r^{2}fl^{2}K^{2}\beta-rfl^{2}K^{2}b\nu-2r^{2}flK\mu\beta\\
-rf\mu^{2}lKb+rf\mu lKb\nu+r^{2}f\mu^{2}\beta+\tau l^{3}K^{3}b^{2}-\tau l^{2}K^{2}b^{2}\mu-\tau l^{2}K^{2}br\beta\\
+\tau l^{2}K^{2}b^{2}\nu+\tau lKb\mu r\beta-l^{3}K^{3}br\beta+l^{3}K^{3}b^{2}\nu+2l^{2}K^{2}b\mu r \beta-l^{2}K^{2}b^{2}\mu \nu\\
-lKb\mu^{2}r\beta),\\
a_2=\frac{1}{l^{2}K^{2}b^{2}}(r^{2}fl{2}K^{2}\beta
-rfl^{2}K^{2}b\nu-2r^{2}flK\mu\beta-r^{2}flK\mu b+rf\mu lKb\nu\\
+r^{2}f\mu^{2}\beta+r^{2}f\mu^{2}b-\tau l^{2}K^{2}br\beta
+\tau l^{2}K^{2}b^{2}\nu+\tau lKb\mu r\beta+\tau lKb^{2}\mu r\\
+\mu r b^{2}l^{2}K^{2}-\mu^{2}rb^{2}lK-\mu r^{2}blK\beta+\mu rb^{2}lK\nu
+\mu^{2}r^{2}b\beta),\\
a_1=\frac 1{lkb}(-rflK+rf\mu+\tau lKb-lKr\beta+lKbv+\mu r \beta).
\end{eqnarray*}
In this case the analysis of the Routh-Hurwitz conditions is far from being easy.

\begin{figure}[ht]
\centering
\includegraphics[width=10cm]{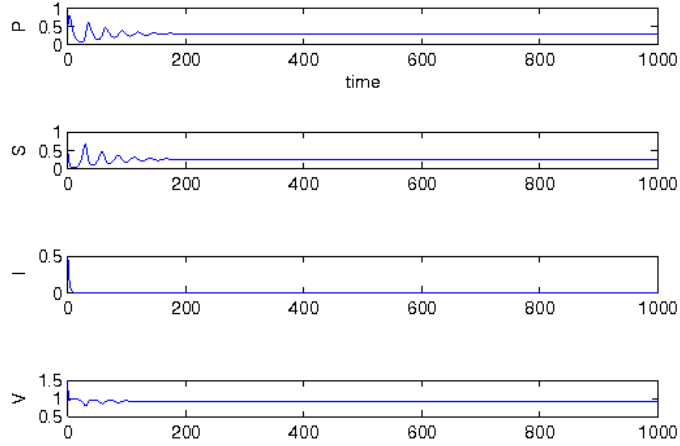}
\caption{The disease-free equilibrium $E_4=(0.3, 0.25, 0, 1.0)$, attained in the logistic case
for the parameter values
$g=0.3$, $f=0.2$, $c=0.4$, $l=0.6$, $q=0.7$, $b=0.9$, $\beta=0.1$, $\tau=0.2$, $\nu=0.2$,
$\mu=0.2$, $r=1.3$, $K=1.0$.}
\label{fig:2}
\end{figure}

\begin{figure}[ht]
\centering
\includegraphics[width=10cm]{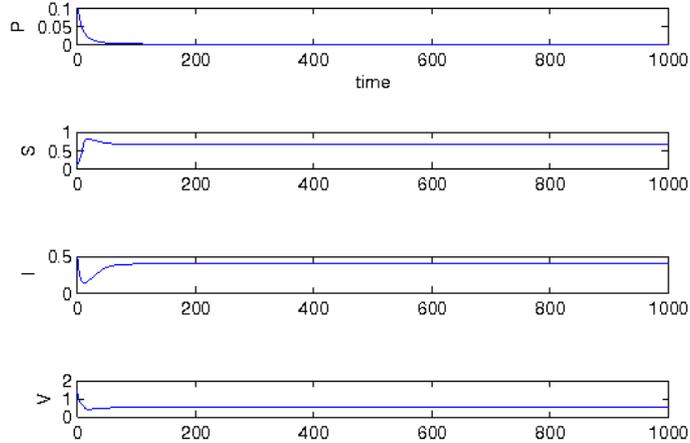}
\caption{The top predator-free equilibrium $E_3=(0,0.6667,0.4103,0.5385)$
stably attained for the parameter values
$g=0.3$, $f=0.1$, $c=0.4$, $l=0.6$, $q=0.7$, $b=0.9$, $\beta=0.3$, $\tau=0.2$, $\nu=0.2$,
$\mu=0.2$, $r=1.3$, $K=1.0$.
}
\label{fig:e12}
\end{figure}

\begin{figure}[ht]
\centering
\includegraphics[width=10cm]{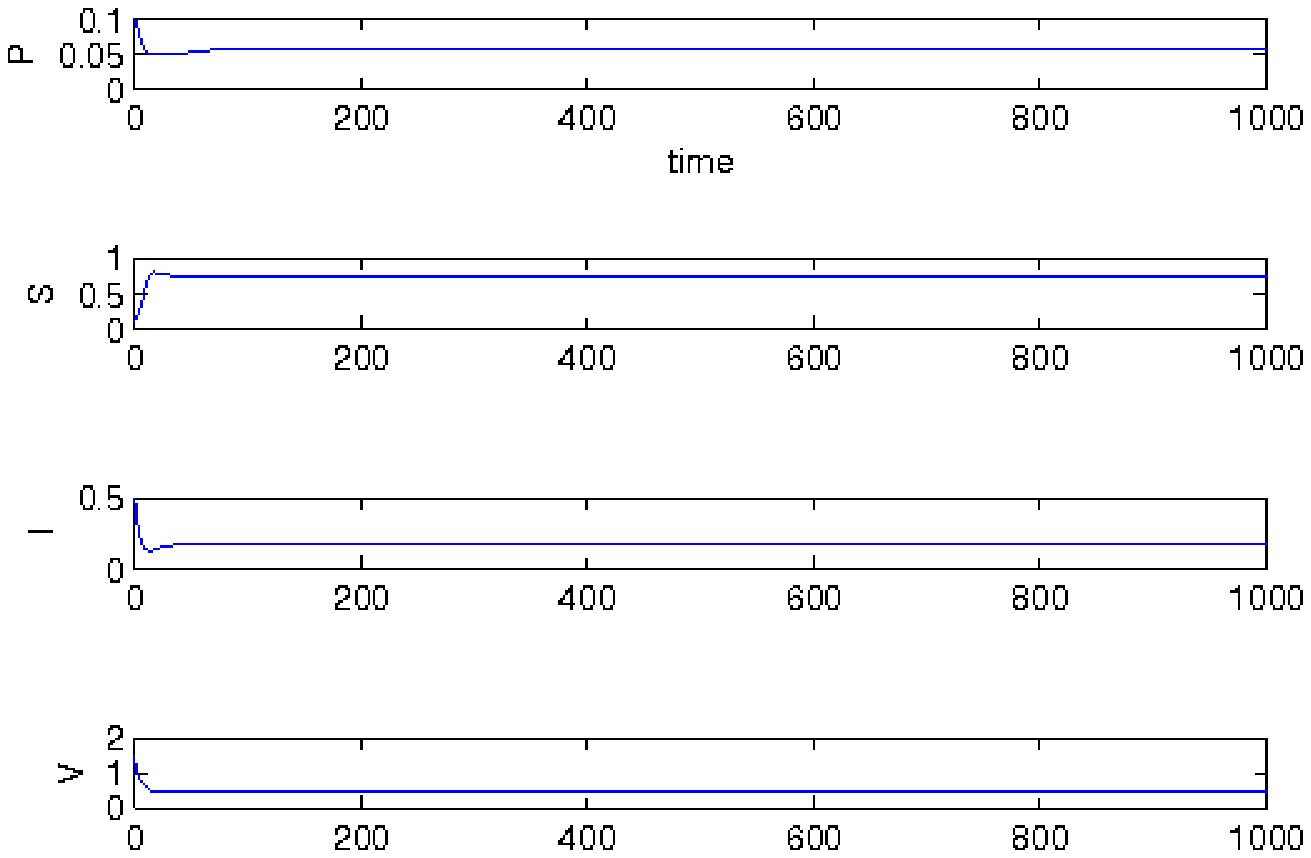}
\caption{The coexistence equilibrium $E^*= (0.0571, 0.7429, 0.1714, 0.4857)$
stably attained in the logistic model for the parameter values
$g=0.3$, $f=0.2$, $c=0.4$, $l=0.6$, $q=0.7$, $b=0.9$, $\beta=0.3$, $\tau=0.2$, $\nu=0.2$,
$\mu=0.2$, $r=1.3$, $K=1.0$.
}
\label{fig:e16}
\end{figure}

The case of the equilibria $E_{3}$ and $E_{4}$ leads to similar very complicated expressions for the coefficients of the characteristic
polynomial (\ref{char_pol}), which we omit altogether.
In the simulations, we show that all these last three points can be attained by the system's trajectories
at a stable level, for suitable parameter choices. These are illustrated in Figures
\ref{fig:2}, \ref{fig:e12}, \ref{fig:e16}.

\section {Discussion}

We have investigated a three trophic level ecoepidemic food chain, in which the disease
affects the population at the intermediate trophic level. In all models the origin is
unstable, this essentially stems from the demographic assumptions, and represents a good
property of the ecosystem, showing that it cannot be completely wiped out.

The purely demographic model admits the following equilibria:
the bottom prey-only equilibrium, which however exists at a finite level only in
the logistic case, the equilibrium with the bottom prey and the
intermediate predator, and coexistence.
These equilibria are related to each
other via two transcritical bifurcations, which occur whenever the parameters
$\rho_1$ and $\rho_2$ cross the critical value 1. In those cases, the intermediate predator
and top predator respectively enter permanently into the system.

In the ecoepidemic models again these purely demographic, disease-free,
equilibria can be found, in particular we observe again that the bottom
prey-only equilibrium exists just in the logistic version.

A very interesting
situation occurs in (\ref{sistMalthus}).
The demographic coexistence equilibrium in the Malthus version
of the ecoepidemic model is not found. The latter is always unstable, so that
the three populations persist with an endemic disease only via sustained oscillations.
Thus introducing a transmissible disease in
a food chain model of this type has the effect that the disease either enters
endemically in the system, or it removes
one trophic level, specifically the uppermost one, if the stability conditions of
equilibrium $\widetilde E_1$ are satisfied, namely (\ref{stab_E7M}).

Alternatively, we can rephrase this concept in a different way.
The only possibility for the disease
to be endemic occurs whenever the coexistence equilibrium is attained. No subsystem
allows the disease to be present in it. Therefore in this system
the task of eradicating the epidemics
is intimately tied to the disappearance of at least one trophic level. Specifically,
it will be the top predator, at the stable equilibrium $\widetilde E_1$,
or possibly at the closed orbits centered around it. This is counterintuitive, since one
expects the top predator to have a positive role in the disease containment.
In fact naively we could think that by
hunting diseased individuals it would contain the epidemics spread.

The same result instead does not hold for the logistic system, we find indeed both the
three-level disease-free food chain equilibrium $E_4$ and the subsystem made of the
lowest two trophic levels with endemic disease, equilibrium $E_3$.

Upon comparison of the two situations, we can conclude thus that
providing more food for the bottom prey, i.e. driving the logistic system toward
the Malthus model, may help in disease eradication, but also may drive to extinction the
top predator. This could be regarded as an alternative formulation of the paradox
of enrichment, by which by feeding the prey one kills the predators.
However, this phenomenon in the present situation is to be ascribed to the demographic
model and not to the disease, as the same occurs in the epidemic-free model. In it,
by providing large amount of food for the bottom prey the bottom prey-only equilibrium $D_1$
disappears, and the system can settle either to coexistence or to the top predator-free
equilibrium $D_2$, depending on the value of the critical parameter $\rho_2$, since the stability
conditions (\ref{D*_stab}) are always satisfied for the coexistence equilibrium. Instead, in
the ecoepidemic model either all the populations of the system, infected included, oscillate,
or the top predator-free equilibrium is attained by imposing condition 
(\ref{Dhat_feas}).

\end{document}